\documentclass[12pt,reqno]{amsart}
\usepackage{amssymb}
\usepackage{amsxtra}

\newcommand{\barF}{\overline{F}}
\newcommand{\e}{\varepsilon}
\newcommand{\G}{\Gamma}
\newcommand{\al}{\alpha}
\newcommand{\eop}{\hfill$\square$}

\theoremstyle{plain}
\newtheorem{Thm}{Theorem}

\newtheorem{Lem}{Lemma}
\theoremstyle{definition}

\theoremstyle{remark}

\newtheorem*{Notn}{Notation}
\newtheorem{Ex}{Example}

\numberwithin{equation}{section}

\begin{document}


\title[Representing the Mean Residual Life]{On Representing the
Mean Residual Life in Terms of the Failure Rate}

\date{\today}

\author{Ramesh~C. Gupta}
\address{Department of Mathematics \& Statistics\\
         University of Maine\\
         5752 Neville Hall
         Orono, Maine 04469-5752\\
         U.S.A.}
\author{David~M. Bradley}
\email[Ramesh~C. Gupta]{rcgupta@maine.maine.edu}

\email[David~M. Bradley]{bradley@math.umaine.edu,
dbradley@member.ams.org}

\subjclass{Primary: 62N02; Secondary: 62N05}

\keywords{Pearson family of distributions, Gaussian probability
function, ultimately increasing failure rate.}

\begin{abstract}
  In survival or reliability studies, the mean residual life or
  life expectancy is an important characteristic of the model.
  Whereas the failure rate can be expressed quite simply in terms
  of the mean residual life and its derivative, the inverse
  problem---namely that of expressing the mean residual life in
  terms of the failure rate---typically involves an integral of a
  complicated expression.  In this paper, we obtain simple expressions
  for the mean residual life in terms of the failure rate for
  certain classes of distributions which subsume many of the
  standard cases.  Several results in the literature can be obtained
  using our approach.  Additionally, we develop an expansion for
  the mean residual life in terms of Gaussian probability
  functions for a broad class of ultimately increasing failure rate
  distributions.  Some examples are provided to illustrate the
  procedure.
\end{abstract}

\maketitle

\section{Introduction}\label{sect:Intro}
In life testing situations, the expected additional lifetime given
that a component has survived until time $t$ is a function of $t$,
called the mean residual life.  More specifically, if the random
variable $X$ represents the life of a component, then the mean
residual life is given by $m(t)=E(X-t|X>t)$.  The mean residual
life has been employed in life length studies by various authors,
e.g.~Bryson and Siddiqui (1969), Hollander and Proschan (1975),
and Muth (1977).  Limiting properties of the mean residual life
have been studied by Meilijson (1972), Balkema and de Hann (1974),
and more recently by Bradley and Gupta (2002).  A smooth estimator
of the mean residual life is given by Chaubey and Sen (1999).

It is well known that the failure rate can be expressed quite
simply in terms of the mean residual life and its derivative:
see~\eqref{DE1} below.  However, the inverse problem---namely that
of expressing the mean residual life in terms of the failure
rate---typically involves an integral of a complicated expression.
In this paper, we obtain a simple expression for the mean residual
life in terms of the failure rate for certain classes of
distributions.  Many of the standard cases are subsumed, and
several results in the literature can be obtained using our
approach.  However, the emphasis here is to express the mean
residual life in terms of the failure rate.  For the class of
ultimately increasing failure rate distributions, we also provide
sufficient conditions under which the mean residual life can be
expanded in terms of Gaussian probability functions. Finally, some
examples are presented to illustrate the procedure.

\section{Background and Definitions}\label{sect:Back}

Let $F:[0,\infty)\to[0,\infty)$ be a nondecreasing,
right continuous function with $F(0)=0$,
$\lim_{x\to\infty}F(x)=1$, and let $\nu$ denote the induced
Lebesgue-Stieljes measure. (Equivalently, let $\nu$ be a
probability measure on $[0,\infty)$ and let $F$ be the cumulative
distribution function of $\nu$.) If $X$ is a nonnegative random
variable representing the life of a component having distribution
function $F$, the mean residual life is defined by
\[
   m(t) = E(X-t|X>t) = \frac{1}{\barF(t)}\int_t^\infty
   (x-t)\,d\nu(x),\qquad t\ge 0,
\]
where $\barF=1-F$ is the so-called survival function.  Writing
$x-t=\int_t^x du$ and employing Tonelli's theorem yields the
equivalent formula
\begin{equation}
   m(t)=\frac{1}{\barF(t)}\int_t^\infty\int_t^x du \,d\nu(x)
   =\frac{1}{\barF(t)}\int_t^\infty\int_u^\infty d\nu(x)\, du
   =\frac{1}{\barF(t)}\int_t^\infty \barF(u)\,du,
\label{MRLbarF}
\end{equation}
which is sometimes also used as a definition.  The cumulative
hazard function may be defined by $R=-\log \barF$.
Then~\eqref{MRLbarF} implies that
\begin{equation}
   m(t)
   = \int_0^\infty \exp\big\{R(t)-R(t+x)\big\}\,dx.
   \label{MRLR}
\end{equation}
If $F$ (equivalently, $\nu$) is also absolutely continuous, then
the probability density function $f$ and the failure rate (hazard
function) $r$ are defined almost everywhere by $f=F'$ and
$r=f/\,\barF=R'$, respectively, and then
\begin{equation}
   R(t)
       = -\log\barF(t)
       =-\int_0^t \frac{d\barF(x)}{\barF(x)}
       = \int_0^t r(x)\,dx.
\label{R}
\end{equation}
In view of~\eqref{MRLR} and~\eqref{R}, we have expressed $m$ in
terms of $r$, albeit somewhat indirectly.

Ideally, we'd like to express the mean residual life in terms of
known functions of the failure rate and its derivatives without
the use of integrals. In any case, it is useful to have
alternative representations of the mean residual life. We note
that the converse problem, that of expressing the failure rate in
terms of the mean residual life and its derivatives is trivial,
for~\eqref{MRLbarF} and~\eqref{R} imply that
\begin{equation}
   m'(t)=r(t)m(t)-1.
\label{DE1}
\end{equation}

\section{A General Family of Distributions, Including the Pearson
Family}\label{sect:Pearson}

Consider the family of distributions whose probability density
function $f$ is differentiable.  Write
\begin{equation}
   \frac{f'(x)}{f(x)}=\frac{\mu-x}{g(x)} - \frac{g'(x)}{g(x)},
\label{family}
\end{equation}
where $\mu$ is a constant, and $g$ satisfies the first order
linear differential equation
\[
   g'(x)+\frac{f'(x)}{f(x)}g(x)=\mu -x.
\]
Alternatively, we may view $g$ as given; then $f$ is uniquely
determined by~\eqref{family}.  Clearing the fractions
in~\eqref{family} and integrating yields
\[
   \int_t^\infty x f(x)\,dx = \mu\barF(t)+g(t)f(t),\qquad
   t\ge 0.
\]
In other words,
\begin{equation}
   E(X|X>t) = \mu + g(t) r(t),\qquad t\ge 0,
\label{Efamily}
\end{equation}
or equivalently,
\begin{equation}
   m(t) = \mu-t+g(t)r(t),\qquad t\ge 0.
\label{mfamily}
\end{equation}
Thus, the mean residual life has been expressed in terms of the
failure rate $r$, the given function $g$, and the constant $\mu$.
By appropriately specializing $g$ in~\eqref{mfamily}, one can
obtain many of the important cases that have appeared in the
literature.  We note that a result similar to~\eqref{mfamily} was
obtained by Ruiz and Navarro~(1994), but their emphasis was
different.

\subsection{The Pearson Family}
The quadratic function $g(x)=a_0+a_1x+a_2x^2$ with $a_2\ne-1/2$
yields the Pearson family, special cases of which include the beta
distributions, the gamma distribution, and the normal
distribution. From~\eqref{family}, we have
\begin{equation}
   \frac{f'(x)}{f(x)}=
   \frac{\mu-x}{a_0+a_1x+a_2x^2}-\frac{a_1+2a_2x}{a_0+a_1x+a_2x^2}
   = -\frac{x+d}{A_0+A_1x+A_2x^2},
\label{Pearson}
\end{equation}
say, where $A_j=a_j/(1+2a_2)$ for $j=0,1,2$ and
$d=(a_1-\mu)/(1+2a_2).$  For the Pearson family,~\eqref{Efamily}
gives
\begin{equation}
   E(X|X>t)=\frac{A_1-d}{1-2A_2}+\frac{A_0+A_1t+A_2t^2}{1-2A_2}r(t).
\label{Pfamily}
\end{equation}
See Nair and Sankaran~(1991) and Ruiz and Navarro~(1994).  Note
that since $A_2=a_2/(1+2a_2)$, it is impossible for $1-2A_2$ to
vanish in~\eqref{Pfamily}.

\subsubsection{The Beta Distributions}  With $g(x)=x(1-x)/(a+b)$, we have
the beta distribution of the first kind:
\[
   f(x)= \frac{\G(a+b)}{\G(a)\G(b)}
   {x^{a-1}(1-x)^{b-1}},
   \quad 0<x<1,\quad a>0,\quad b>0.
\]
It satisfies
\[
   \frac{f'(x)}{f(x)}= -\frac{x+(1-a)/(a+b-2)}{x(1-x)/(a+b-2)},
   \qquad a+b\ne 2,
\]
and hence belongs to the Pearson family given by~\eqref{Pearson}
with $d=(1-a)/(a+b-2)$, $A_0=0$, $A_1=-A_2=1/(a+b-2)$, and $a+b\ne
2$. Thus,
\[
   E(X|X>t) = \mu +\frac{t(1-t)}{a+b}r(t), \qquad 0<t<1,
\]
where $\mu=a/(a+b)$.  If $a+b=2$, we have
\[
   f(x)= \frac{1}{\G(a)\G(2-a)}\left(\frac{x}{1-x}\right)^{a-1},
   \qquad 0<a<2,
\]
and $g(x)=x(1-x)/2$, $\mu=a/2$, from which~\eqref{Efamily} gives
\[
   E(X|X>t)= \tfrac12a+\tfrac12t(1-t)r(t),\qquad 0<t<1.
\]

The beta distribution of the second kind has the form
\[
   f(x) = \frac{c x^{\beta-1}}{(\gamma+x)^{\al+\beta}},
   \qquad c,\al,\beta,\gamma,x>0.
\]
Since then
\[
   \frac{f'(x)}{f(x)}= -\frac{x+(1-\beta)\gamma/(\al+1)}
   {x^2/(\al+1)+\gamma x/(\al+1)},
\]
$f$ belongs to the Pearson family~\eqref{Pearson} with
$d=(1-\beta)\gamma/(1+\al)$, $A_0=0$, $A_1=\gamma/(1+\al)$, and
$A_2=1/(1+\al)$. Hence,
\[
   E(X|X>t)=\mu+\left(\frac{t^2+\gamma t}{\al-1}\right)r(t),
   \qquad \al\ne1,
\]
where $\mu=\beta\gamma/(\al-1)$ is the mean. A similar result was
obtained by Ahmed~(1991) using a completely different approach.

\subsubsection{The Gamma Distribution}
Let $B>0$.  The linear function $g(x)=Bx$ for $x>0$ yields the
gamma distribution.  In this case,~\eqref{family} gives
\[
   \frac{f'(x)}{f(x)}=-\frac{x+B-\mu}{Bx},
\]
so that $f$ belongs to the Pearson family~\eqref{Pearson} with
$d=B-\mu$, $A_0=0$, $A_1=B$, and $A_2=0$.  In fact,
\[
   f(x)=\frac{B^{-\mu/B}}{\G(\mu/B)} x^{\mu/B-1}e^{-x/B},
\]
a gamma distribution with mean $\mu$.  Hence~\eqref{Efamily} takes
the form
\[
   E(X|X>t)=\mu+Btr(t)
\]
in this case.  This is Theorem 1 of Osaki and Li~(1988).  See also
El-Arishy~(1995).

\subsubsection{The Normal Distribution}
Let $\sigma>0$.  The constant function $g(x)=\sigma^2$ for
$-\infty<x<\infty$ yields the normal distribution.  In this
case,~\eqref{family} gives
\[
   \frac{f'(x)}{f(x)}=\frac{\mu-x}{\sigma^2},
\]
so that $f$ belongs to the Pearson family~\eqref{Pearson} with
$d=-\mu$, $A_0=\sigma^2$, and $A_1=A_2=0$.  It follows that
\[
   f(x)=\sigma^{-1}(2\pi)^{-1/2}\exp\big\{-(x-\mu)^2/2\sigma^2\big\},
\]
a normal distribution with mean $\mu$ and variance $\sigma^2$.
Hence~\eqref{Efamily} can be written as
\[
   E(X|X>t) = \mu+\sigma^2 r(t).
\]
This is Theorem 3.1 of Ahmed and Abdul-Rahman~(1993).  See also
McGill~(1992) and Kotz and Shanbhag~(1980).

\subsection{The Maxwell Distribution}
The Maxwell distribution has the form
\[
    f(x)=4b^{-3}\pi^{-1/2}x^2\exp\{-x^2/b^2\},
   \qquad x>0,\qquad b>0.
\]
See El-Arishy~(1993).  We then have
\[
   \frac{f'(x)}{f(x)} = \frac{2}{x}-\frac{2x}{b^2},
\]
which can be written in the form~\eqref{family} with $\mu=0$ and
$g(x)=(1+b^2/x^2)b^2/2$.  Hence,
\[
   E(X|X>t) = \frac{(t^2+b^2)b^2 r(t)}{2t^2}.
\]
The fact that the corresponding failure rate is increasing can be
seen by using Glaser's (1980) result.

\section{Ultimately Increasing Failure Rate Distributions}
Consider the class 
of distributions whose failure rate is ultimately increasing. More
specifically, the failure rate should be strictly increasing from
some point onward. Obviously, the important class of lifetime
distributions having a bathtub-shaped failure rate with change
points $0\le t_1\le t_2<\infty$ (i.e.\ for which the failure rate
is strictly decreasing on the interval $[0,t_1]$, constant on
$[t_1,t_2]$ and strictly increasing on $[t_2,\infty)$) constitutes
a proper subclass of the distributions we consider here.
We'll see that if the failure rate is strictly increasing from
some point onward, then under certain additional conditions the
mean residual life can be expanded in terms of Gaussian
probability functions.

\begin{Notn} Our conventions regarding the
Bachmann-Landau $O$-notation, the Vinogradov $\ll$-notation, and
$o$-notation are fairly standard.  Thus, if $h$ is a function of a
positive real variable, the symbol $O(h(t))$, $t\to\infty$,
denotes an unspecified function $g$ for which there exist positive
real numbers $t_0$ and $B$ such that $|g(t)|\le B|h(t)|$ for all
real $t> t_0$.  For such $g$ we write $g(t)\ll h(t)$ or
$g(t)=O(h(t))$.  The notation $g(t)=o(h(t))$, $t\to\infty$, means
that for every real $\e>0$, no matter how small, there exists a
positive real number $t_0$ such that $|g(t)|\le \e |h(t)|$
whenever $t>t_0$.
\end{Notn}
%

\begin{Thm}\label{Thm:normal}
Suppose that from some point onward, the failure rate $r$
increases (strictly) without bound.  Suppose further that for some
positive integer $n$, the $n-1$ derivative is continuous and
satisfies
\begin{equation}
   \big|r^{(n-1)}(t+x)\big| \ll \big|r^{(n-1)}(t)\big|,
   \qquad t\to\infty,
\label{uniform}
\end{equation}
uniformly in $x$ for $0\le x\le \min\big(1, |r''(t)|^{-1/3}\big)$,
and
\[
   r^{(j)}(t) \ll \max\left(1,|r''(t)|^{(j+1)/3}\right),
   \qquad   t\to\infty, \qquad 3\le j\le n-1.
\]
Finally, suppose there exists a positive real number $\e$ such
that for each integer $j$ in the range $3\le j\le n$,
\begin{equation}
   r^{(j-1)}(t) = o\left((r(t))^{j-j\e}\right),\qquad t\to\infty.
\label{eps}
\end{equation}
Then we have the following expansion for the mean residual life:
\begin{equation}
   m(t) = \sum_{k=0}^{n-1} b_k(t)\varphi_k(t)
    +o\left((r(t))^{-1-n\e}\right),
   \qquad  t\to\infty,
\label{normal}
\end{equation}
where the coefficients $b_k(t)$ are given by the formal power
series identity
\begin{equation}
   \sum_{k=0}^\infty b_k x^k = \exp\bigg\{-\sum_{k=3}^\infty
   r^{(k-1)}(t)\frac{x^k}{k!}\bigg\},
\label{bkdef}
\end{equation}
and
\begin{align}\label{phikdef}
\begin{split}
   \varphi_k(t)
   &=\int_0^\infty x^k \exp\{-xr(t)-\tfrac12 x^2r'(t)\}\,dx\\
   &= (-1)^k\sqrt{\frac{2\pi}{r'(t)}}\bigg[\frac{\partial^k}{\partial p^k}
   \bigg(1-\Phi\bigg(\frac{p}{\sqrt{r'(t)}}\bigg)\bigg)
   \exp\bigg\{ \frac{p^2}{2r'(t)}\bigg\}\bigg|_{p=r(t)}\bigg].
\end{split}
\end{align}
Here,
\begin{equation}
   \Phi(x)=\frac{1}{\sqrt{2\pi}}\int_{-\infty}^x
   e^{-v^2/2}\,dv
\label{Phi}
\end{equation}
is the Gaussian probability function, i.e.\ the cumulative
distribution function of the standard normal distribution.
\end{Thm}

Before proving Theorem~\ref{Thm:normal}, we make some preliminary
remarks and give two illustrative examples.  First, if $r''(t)=0$,
then the uniformity condition on $x$ in~\eqref{uniform} should be
interpreted as $0\le x\le 1.$  Next, observe that the
hypothesis~\eqref{eps} becomes more restrictive as $\e$ increases.
In particular $\e>1$ implies $\lim_{t\to\infty}r^{(j-1)}(t)=0$. Of
course, as $\e$ increases, the error term in~\eqref{normal}
decreases.  On the other hand, if $0<\e<1$, then $r^{(j-1)}(t)$
does not necessarily approach zero, but the correspondingly weaker
hypothesis implies a weaker conclusion (larger error term). In any
case, since $r$ increases without bound, the error term tends to
zero as $t\to\infty$. Additionally, if $r$ is infinitely
differentiable we may let $n\to\infty$ in~\eqref{normal} to obtain
the \emph{convergent} infinite series expansion
\[
    m(t) = \sum_{k=0}^{\infty} b_k(t)\varphi_k(t),
\]
valid for all sufficiently large values of $t$.  (More
specifically, for those $t$ for which $r(t)>1$.)

In general, however, we do not assume the failure rate has
infinitely many derivatives; $n$ is fixed and the generating
function~\eqref{bkdef} is a \emph{formal} power series.
Expanding~\eqref{bkdef} to compute $b_k(t)$ in terms of $r(t)$ and
its derivatives shows that if $r$ has only $n-1$ derivatives, then
$b_k(t)$ is undefined if $k\ge n$. Differentiating~\eqref{bkdef}
leads to the recurrence
\begin{equation}
   b_{k+1}(t) = -\frac{1}{k+1}\sum_{j=2}^k
   \frac{r^{(j)}(t)}{j!}b_{k-j}(t),\qquad k\ge 2,
\label{recurr}
\end{equation}
from which the coefficients $b_k(t)$ may be successively
determined, starting with the initial values $b_0(t)=1$,
$b_1(t)=b_2(t)=0$. On the other hand, an application of the
multinomial theorem yields the explicit representation
\begin{equation}
   b_k(t) = \sum_{p=0}^{\lfloor k/3 \rfloor} (-1)^p
         \sum \prod_{j\ge 2}
         \frac{1}{\al_j!}\bigg(\frac{r^{(j)}(t)}{(j+1)!}\bigg)^{\al_j},
\label{explicit}
\end{equation}
in which $\lfloor k/3\rfloor$ is the greatest integer not
exceeding $k/3$ and the inner sum is over all non-negative
integers $\al_2,\al_3,\dots$ such that $\sum_{j\ge 2}\al_j=p$ and
$\sum_{j\ge 2}(j+1)\al_j=k$.

Finally, we note that the functions $\varphi_k$ of~\eqref{phikdef}
may also be given more explicitly.  By setting $a=r(t)$ and
$b=2r'(t)$ in Lemma 1 below, we find that
\begin{multline*}
 \varphi_k(t)
   = (-1)^k \bigg(\frac{2}{r'(t)}\bigg)^{(k+1)/2}\\
   \times\bigg\{
   \sum_{h=0}^{\lfloor k/2\rfloor}
   \binom{k}{2h}\lambda^{k/2-h}\,\G(h+1/2)
   \bigg(\big(1-\Phi\big(\sqrt{2\lambda}\big)\big)e^{\lambda}+
   \frac12\sum_{j=0}^{h-1}\frac{\lambda^{j+1/2}}{\G(j+3/2)}
   \bigg)\\
   -\frac12
   \sum_{h=0}^{\lfloor
   k/2\rfloor}h!\binom{k}{2h+1}\lambda^{(k-1)/2-h}\sum_{j=0}^h
   \frac{\lambda^j}{j!}\bigg\},
\end{multline*}
where $\lambda= (r(t))^2/2r'(t)$, $\G(h+1/2) =
\pi^{1/2}\prod_{j=1}^h (j-1/2)$, and $\Phi$ denotes the Gaussian
probability function~\textup{\eqref{Phi}}.

\begin{Lem}\label{Lem:normal} Let $a$ be a real number, let $b$ be a
positive real number, and let $k$ be a non-negative integer.  Then
\begin{multline*}
   \int_0^\infty x^k\exp\{-ax-bx^2\}\,dx\\
   = (-1)^k b^{-(k+1)/2}  \bigg\{
   \sum_{h=0}^{\lfloor k/2\rfloor}
   \binom{k}{2h}\lambda^{k/2-h}\,\G(h+1/2)
     \bigg(\big(1-\Phi\big(\sqrt{2\lambda}\big)\big)e^{\lambda}+
   \frac12\sum_{j=0}^{h-1}\frac{\lambda^{j+1/2}}{\G(j+3/2)}
   \bigg)\\
   -\frac12
   \sum_{h=0}^{\lfloor
   k/2\rfloor}h!\binom{k}{2h+1}\lambda^{(k-1)/2-h}\sum_{j=0}^h
   \frac{\lambda^j}{j!}\bigg\},
\end{multline*}
where $\lambda= a^2/4b$.
\end{Lem}
Proofs of Theorem~\ref{Thm:normal} and Lemma~\ref{Lem:normal} are
relegated to \S\ref{section:proofs}, the final section.

\section{Applications}

We provide two examples indicating how Theorem~\ref{Thm:normal}
may be applied.

\begin{Ex}  Consider a linear failure rate of the form
\begin{equation}
   r(t)=\al+\beta t,\qquad \beta>0.
\label{LFR}
\end{equation}
The motivation and application of~\eqref{LFR} to analyzing various
data sets has been demonstrated by Kodlin (1967) and Carbone et
al.~(1967).  Statistical inference related to the linear failure
rate model has been studied by Bain (1974), Shaked (1974) and more
recently by Sen and Bhatacharya (1995).  For this model, the
hypotheses of Theorem~\ref{Thm:normal} are trivially satisfied for
\emph{any} positive integer $n$ and \emph{any} positive real
number $\e$.  Since $r''$ vanishes identically in this case, we
see that $b_k(t)=0$ for $k>0$ in~\eqref{normal} and in fact we
have the exact result
\[
   m(t)=\int_0^\infty \exp\{-(\al+\beta t)x-\beta x^2/2\}\,dx
   = \exp\bigg\{\frac{(\al+\beta t)^2}{2\beta}\bigg\}
   \bigg(1-\Phi\bigg(\frac{\al+\beta
   t}{\sqrt{\beta}}\bigg)\bigg)\sqrt{\frac{2\pi}{\beta}}.
\]
\end{Ex}

\begin{Ex}  Chen~(2000) proposes the two-parameter distribution
with cumulative distribution function given by
\[
   F(t) = 1 - \exp\big\{\big(1-\exp(t^\beta)\big)\lambda\big\},
   \qquad t>0,
\]
where $\lambda>0$ and $\beta>0$ are parameters.  The corresponding
hazard function is the ultimately strictly increasing function of
$t$ given by
\begin{equation}
   r(t) = \lambda\beta t^{\beta-1}\exp(t^\beta),
   \qquad t>0.
\label{ChenHazard}
\end{equation}
It is straightforward, albeit somewhat tedious, to verify that
Chen's failure rate~\eqref{ChenHazard} satisfies the hypotheses of
Theorem~\ref{Thm:normal} with $n>2$ and $0<\e\le 2/3$. Clearly
$\e=2/3$ is optimal here. Thus, with derivatives of $r$
in~\eqref{bkdef} and~\eqref{phikdef} now coming
from~\eqref{ChenHazard}, we see that the asymptotic formula
\[
    m(t) = \sum_{k=0}^{n-1} b_k\varphi_k(t)
    +o\left((r(t))^{-1-2n/3}\right),
   \qquad  t\to\infty,
\]
holds for \emph{all} integers $n>2$. In particular, as the error
term in~\eqref{normal} tends to zero in the limit as $n\to\infty$,
we obtain the \emph{convergent} infinite series representation
\[
   m(t) = \sum_{k=0}^\infty b_k(t)\varphi_k(t),
\]
valid for all sufficiently large values of $t$.  There is no need
to work out the coefficients $b_k(t)$ explicitly in this case. One
can simply use the recurrence~\eqref{recurr} to generate them.
\end{Ex}

\section{Proofs}\label{section:proofs}
\subsection{Proof of Theorem~\ref{Thm:normal}}
Since the failure rate is strictly increasing from some point
onward, there exists $t_0\ge 0$ such that $r'(t)>0$ for all $t\ge
t_0$. Also, since $\lim_{t\to\infty}r(t)=\infty$, there exists
$t_1\ge 0$ such that $r(t)\ge 1$ for $t\ge t_1$.  Now let $t\ge
\max(t_0,t_1)$,
$\delta=\delta(t)=\min\big(1,1/\sqrt[3]{|r''(t)|}\big)$, and set
\[
   I(t) := \int_0^\delta \exp\{R(t)-R(t+x)\}\,dx,
   \qquad J(t) := \int_{\delta}^\infty \exp\{R(t)-R(t+x)\}\,dx,
\]
so that $m(t)=I(t)+J(t).$  We have
\begin{align*}
   J(t)
   &\le \int_{\delta}^\infty r(t+x)\exp\left\{R(t)-R(t+x)\right\}\,dx\\
   &= \exp\left\{R(t)-R(t+\delta)\right\}\\
   &= \exp\bigg\{-\delta r(t)-\int_0^\delta x
   r'(t+\delta-x)\,dx\bigg\}\\
   &\le \exp\left\{-\delta r(t)\right\}.
\end{align*}
But $r''=o(r^{3-3\e})$.  By definition of $\delta$, it follows
that from some point onward we must have $\delta r \ge \min (r,
r^{\e})$. Therefore, if we set $\nu = \min(1,\e)$, then $\nu>0$
and
\begin{equation}
   J(t) \le \exp\big\{-(r(t))^\nu\big\},
\label{Jest}
\end{equation}
for all sufficiently large values of $t$.

Next, we write
\begin{align*}
   I(t) &= \int_0^\delta \exp\bigg\{-\sum_{k=1}^{n-1} r^{(k-1)}(t)
   \frac{x^k}{k!}-\frac{1}{(n-1)!}\int_0^x u^{n-1}
   r^{(n-1)}(t+x-u)\,du\bigg\}\,dx\\
   &= \int_0^\delta \bigg(\sum_{k=0}^{n-1} b_k(t) x^k
     +E_n(x,t)\bigg)\exp\left\{-xr(t)-\tfrac12
     x^2r'(t)\right\}\,dx,
\end{align*}
where
\begin{align*}
   E_n(x,t) & =
   \exp\left\{-\sum_{k=3}^{n-1}r^{(k-1)}(t)\frac{x^k}{k!}
   -\frac{1}{(n-1)!}\int_0^x
   u^{n-1}r^{(n-1)}(t+x-u)\,du \right\}
   -\sum_{k=0}^{n-1}b_k(t)x^k\\
    &= O\bigg(x^n\max \prod_{j=2}^{n-1} \big|r^{(j)}(t)\big|^{\al_j}
   \bigg),
\end{align*}
and the maximum is taken over all non-negative integers $\al_j$
satisfying $\sum_{j=2}^{n-1} (j+1)\al_j = n.$  In view of the fact
that $r^{(j-1)}=o(r^{j-j\e})$ for $3\le j\le n$, it follows that
\[
   E_n(x,t) = o(x^n (r(t))^{n-n\e}),\qquad 0\le x\le \delta.
\]
If we now write
\begin{multline*}
   I(t) = \sum_{k=0}^{n-1}b_k(t)\int_0^{\infty} x^k
          \exp\left\{-xr(t)-\tfrac12 x^2r'(t)\right\}\,dx\\
        - \sum_{k=0}^{n-1}b_k(t)\int_{\delta}^\infty x^k
          \exp\left\{-xr(t)-\tfrac12 x^2r'(t)\right\}\,dx\\
          + \int_0^{\delta} E_n(x,t)
          \exp\left\{-xr(t)-\tfrac12 x^2r'(t)\right\}\,dx,
\end{multline*}
then we find that
\[
   I = \sum_{k=0}^{n-1} b_k\varphi_k
        + \sum_{k=0}^{n-1} O\big(b_k r^k e^{-\delta r}\big)
        + o\bigg(r^{n-n\e}\int_0^{\delta} x^n
        e^{-xr}\,dx\bigg).
\]
The hypotheses on $r$ and the definition of the coefficients $b_k$
imply that $b_k=O(r^{k-k\e})$. From the derivation of the
estimate~\eqref{Jest} for $J$ we recall that $\exp(-\delta r)\le
\exp(-r^{\nu})$, at least from some point onward. Finally, as
\[
   \int_0^\delta x^n e^{-xr}\,dx \le \int_0^\infty x^n e^{-xr}\,dx
   = n!\, r^{-n-1},
\]
it follows that
\begin{equation}
   I = \sum_{k=0}^{n-1}b_k\varphi_k + o\big(r^{-1-n\e}\big).
\label{Iest}
\end{equation}
Since $m=I+J$, combining~\eqref{Jest} and~\eqref{Iest} gives the
stated result for $m$.

To complete the proof, it remains only to establish the asserted
evaluation of the integrals $\varphi_k$.  But this is readily
obtained by completing the square in the exponential and
differentiating under the integral. \eop

\subsection{Proof of Lemma~\ref{Lem:normal}} By a straightforward change
of variables, we find that
\begin{align}
   2 b^{(k+1)/2}\int_0^\infty x^k \exp\{-ax-bx^2\}\,dx
   &=  e^{\lambda}
   \int_{\lambda}^\infty \big(t^{1/2}-\lambda^{1/2}\big)^k e^{-t}
   t^{-1/2}\,dt\nonumber\\
   &= e^{\lambda}
   \sum_{h=0}^k\binom{k}{h}(-1)^{k-h}\lambda^{(k-h)/2}\,
   \G\big(\tfrac{h+1}{2},\lambda\big),
\label{incomplete}
\end{align}
where
\begin{equation}
   \G(\al,\lambda) = \int_{\lambda}^\infty t^{\al-1}e^{-t}\,dt
\label{IncGamma}
\end{equation}
is the incomplete gamma function.  If we
integrate~\eqref{IncGamma} by parts and then divide both sides of
the result by $\G(\al+1)=\al\G(\al)$, we obtain the recurrence
formula
\[
   \frac{\G(\al+1,\lambda)}{\G(\al+1)}
   = \frac{\lambda^{\al}e^{-\lambda}}{\G(\al)}
   +\frac{\G(\al,\lambda)}{\G(\al)},
\]
which can be iterated to give
\begin{equation}
   \frac{\G(\al+k,\lambda)}{\G(\al+k)}
   =
   e^{-\lambda}\sum_{h=0}^{k-1}\frac{\lambda^{\al+h}}{\G(\al+h+1)}
   +\frac{\G(\al,\lambda)}{\G(\al)},
\label{IncSum}
\end{equation}
valid for any non-negative integer $k$.  In particular, when
$\al=0$,
\begin{equation}
   \frac{\G(k,\lambda)}{\G(k)} = e^{-\lambda}\sum_{h=0}^{k-1}
   \frac{\lambda^h}{h!}.
\label{IncInteger}
\end{equation}
Equation~\eqref{IncInteger} is valid for all positive integers $k$
if $\lambda\ge 0$; it is also valid when $k=0$ if $\lambda>0$.

Substituting $\al=1/2$ in~\eqref{IncSum} yields
\begin{align}
   \frac{\G(k+1/2,\al)}{\G(k+1/2)}
   &=e^{-\lambda}\sum_{h=0}^{k-1}\frac{\lambda^{k+1/2}}{\G(k+3/2)}
   +\frac{\G(1/2,\lambda)}{\G(1/2)}\nonumber\\
   &=e^{-\lambda}\sum_{h=0}^{k-1}\frac{\lambda^{k+1/2}}{\G(k+3/2)}
   + 2\left(1-\Phi\big(\sqrt{2\lambda}\big)\right).
\label{IncHalfInteger}
\end{align}
Using~\eqref{IncInteger} and~\eqref{IncHalfInteger}, we get the
stated result from~\eqref{incomplete}. \eop

\end{document}